\documentclass[11pt]{article}
\usepackage{amssymb}
\usepackage{amsmath}
\usepackage[all]{xy}

\title{\bf $k$-Gorenstein Modules
}
\author{Zhaoyong Huang\thanks{\small \it E-mail address: huangzy@nju.edu.cn}\\
{\small \it Department of Mathematics, Nanjing University,}\\
{\small \it Nanjing 210093, People's Republic of China}\\}
\usepackage{amssymb}
\date{}
\begin{document}
\baselineskip=18pt \maketitle

\begin{abstract}
Let $\Lambda$ and $\Gamma$ be artin algebras and $_{\Lambda}U
_{\Gamma}$ a faithfully balanced selforthogonal bimodule. In this
paper, we first introduce the notion of $k$-Gorenstein modules
with respect to $_{\Lambda}U _{\Gamma}$ and then characterize it
in terms of the $U$-resolution dimension of some special injective
modules and the property of the functors ${\rm Ext}^{i}({\rm
Ext}^{i}(-, U), U)$ preserving monomorphisms, which develops a
classical result of Auslander. As an application, we study the
properties of dual modules relative to Gorenstein bimodules. In
addition, we give some properties of $_{\Lambda}U _{\Gamma}$ with
finite left or right injective dimension.
\end{abstract}

\vspace{0.5cm}

\centerline{\bf 1. Introduction}

\vspace{0.2cm}

Let $\Lambda$ be a ring. We use mod $\Lambda$ (resp. mod $\Lambda
^{op}$) to denote the category of finitely generated left
$\Lambda$-modules (resp. right $\Lambda$-modules).

Let $\Lambda$ and $\Gamma$ be rings. A bimodule $_{\Lambda}T
_{\Gamma}$ is said to be faithfully balanced if the natural maps
$\Lambda \rightarrow {\rm End}(T_{\Gamma})$ and $\Gamma
\rightarrow {\rm End}(_{\Lambda}T)^{op}$ are isomorphisms; and it
is said to be selforthogonal if ${\rm
Ext}_{\Lambda}^{i}(_{\Lambda}T , {_{\Lambda}T})=0$ and ${\rm
Ext}_{\Gamma}^{i}(T_{\Gamma}, T_{\Gamma})=0$ for any $i \geq 1$.

Let $U$ and $A$ be in mod $\Lambda$ (resp. mod $\Gamma ^{op}$) and
$i$ a non-negetive integer. We say that the grade of $A$ with
respect to $U$, written grade$_{U}A$, is greater than or equal to
$i$ if Ext$_{\Lambda}^{j}(A, U)=0$ (resp. Ext$_{\Gamma}^{j}(A, U
)=0$) for any $0 \leq j <i$. We say that the strong grade of $A$
with respect to $U$, written s.grade$_{U}A$, is greater than or
equal to $i$ if grade$_{U}B \geq i$ for all submodules $B$ of $A$
(see [11]). We give the definition of ($k$-)Gorenstein modules in
terms of strong grade of modules as follows.

\vspace{0.2cm}

{\bf Definition 1.1} For a non-negative integer $k$, a module $U
\in$ mod $\Lambda$ with $\Gamma ={\rm End}(_{\Lambda}U)$ is called
$k$-Gorenstein if s.grade$_{U}{\rm Ext}^{i}_{\Gamma} (N, U)\geq i$
for any $N \in$mod $\Gamma ^{op}$ and $1 \leq i \leq k$. $U$ is
called Gorenstein if it is $k$-Gorenstein for all $k \geq 1$.
Similarly, we may define the notions of $k$-Gorenstein modules and
Gorenstein modules in mod $\Gamma ^{op}$. A bimodule $_{\Lambda}U
_{\Gamma}$ is called a ($k$-)Gorenstein bimodule if both
$_{\Lambda}U$ and $U _{\Gamma}$ are ($k$-)Gorenstein.

\vspace{0.2cm}

{\bf Definition 1.2}$^{[18]}$ Let $U$ be in mod $\Lambda$ (resp.
mod $\Gamma ^{op}$) and $k$ a non-negetive integer. A module $M$
in mod $\Lambda$ (resp. mod $\Gamma ^{op}$) is said to have
$U$-dominant dimension greater than or equal to $k$, written
$U$-dom.dim$(_{\Lambda}M)$ (resp. $U$-dom.dim$(M_{\Gamma}))\geq
k$, if each of the first $k$ terms in a minimal injective
resolution of $M$ is cogenerated by $_{\Lambda}U$ (resp. $U
_{\Gamma}$), that is, each of these terms can be embedded into a
direct product of copies of $_{\Lambda}U$ (resp. $U _{\Gamma}$).

\vspace{0.2cm}

It is clear that any module in mod $\Lambda$ (resp. mod $\Gamma
^{op}$) is 0-Gorenstein. Let $\Lambda$ and $\Gamma$ be artin
algebras and $_{\Lambda}U _{\Gamma}$ a faithfully balanced
selforthogonal bimodule with $_{\Lambda}U \in$mod $\Lambda$ and
$U_{\Gamma} \in$mod $\Gamma ^{op}$. If $U$-dom.dim$(_{\Lambda}U)
\geq k$, then each of the first $k$ terms in a minimal injective
resolution of $_{\Lambda}U$ is finitely cogenerated, and so each
of these terms can be embedded into a finite direct product of
copies of $_{\Lambda}U$. It follows from Lemma 2.6 below that
$_{\Lambda}U$ is $k$-Gorenstein. It was showed in [13] that
$U$-dom.dim$(_{\Lambda}U)=U$-dom.dim$(U_{\Gamma})$. So, at this
moment, $U_{\Gamma}$ is also $k$-Gorenstein. Recall from [20] that
a module $M$ in mod $\Lambda$ (resp. mod $\Gamma ^{op}$) is called
a QF-3 module if $G(M)$ has a cogenerator which is a direct
summand of every other cogenerator, where $G(M)$ is the
subcategory of mod $\Lambda$ (resp. mod $\Gamma ^{op}$) consisting
of all submodules of the modules generated by $M$. It was showed
in [20] Proposition 2.2 that a finitely cogenerated
$\Lambda$-module (resp. $\Gamma ^{op}$-module) $M$ is a QF-3
module if and only if $M$ cogenerates its injective envelope. So
we have that $_{\Lambda}U$ (resp. $U _{\Gamma}$) is 1-Gorenstein
if it is a QF-3 module.

A left and right noetherian ring $\Lambda$ is called
$k$-Gorenstein if for any $1 \leq i \leq k$ the $i$th term in a
minimal injective resolution of $_{\Lambda}\Lambda$ has flat
dimension at most $i-1$. This notion was introduced by Auslander
and Reiten in [5] as a non-commutative version of commutative
Gorenstein rings. By Definition 1.1 and [7] Auslander's Theorem
3.7, $\Lambda$ is a $k$-Gorenstein ring if it is $k$-Gorenstein as
a $\Lambda$-module. Auslander further proved that the notion of
$k$-Gorenstein rings is left-right symmetric (see [7] Auslander's
Theorem 3.7). Wakamatsu in [19] Theorem 7.5 generalized this
result and established the left-right symmetry of the notion
$k$-Gorenstein modules.

In this paper, we will give some further characterizations of
$k$-Gorenstein modules in terms of the $U$-resolution dimension of
some special injective modules and the property of the functors
${\rm Ext}^{i}({\rm Ext}^{i}(-, U), U)$ preserving monomorphisms,
which develops the result of Auslander mentioned above. Our
characterizations will lead to a better comprehension about the
theory of selforthogonal bimodules and cotilting theory (note: the
class of cotilting bimodules is such a kind of faithfully balanced
selforthogonal bimodules with finite left and right injective
dimensions$^{[9]}$).

Throughout this paper, both $\Lambda$ and $\Gamma$ are artin
algebras (unless stated otherwise), $_{\Lambda}U _{\Gamma}$ is a
faithfully balanced selforthogonal bimodule with $_{\Lambda}U
\in$mod $\Lambda$ and $U_{\Gamma} \in$mod $\Gamma ^{op}$.

The following is an outline of this paper. In Section 2 we list some
lemmas which will be used later. In Section 3 we characterize
$k$-Gorenstein modules with respect to $_{\Lambda}U _{\Gamma}$ in
terms of the $U$-resolution dimension (see Section 2 for the
definition) of some special injective modules and the property of
the functors ${\rm Ext}^{i}({\rm Ext}^{i}(-, U), U)$ preserving
monomorphisms. In fact, we will prove the following theorem, which
extends [7] Auslander's Theorem 3.7.

{\bf Theorem} {\it The following statements are equivalent.}

(1) $_{\Lambda}U$ {\it is} $k$-{\it Gorenstein.}

(2) $U$-resol.dim$_{\Lambda}(E_{i}) \leq i$, {\it where} $E_i$
{\it is the} $(i+1)$st {\it term in a minimal injective of} $U$
{\it as a left} $\Lambda$-{\it module, for any} $0 \leq i \leq
k-1$.

(3) Ext$_{\Gamma}^{i}({\rm Ext}_{\Lambda}^{i}(-, U), U)$: mod
$\Lambda \to$ mod $\Lambda$ {\it preserves monomorphisms} {\it for
any} $0 \leq i \leq k-1$.

(1)$^{op}$ $U _{\Gamma}$ {\it is} $k$-{\it Gorenstein.}

(2)$^{op}$ $U$-resol.dim$_{\Gamma}(E'_{i}) \leq i$, {\it where}
$E'_i$ {\it is the} $(i+1)$st {\it term in a minimal injective of}
$U$ {\it as a right} $\Gamma$-{\it module, for any} $0 \leq i \leq
k-1$.

(3)$^{op}$ Ext$_{\Lambda}^{i}({\rm Ext}_{\Gamma}^{i}(-, U), U)$:
mod $\Gamma ^{op} \to$ mod $\Gamma ^{op}$ {\it preserves
monomorphisms} {\it for any} $0 \leq i \leq k-1$.

As mentioned above, Wakamatsu in [19] Theorem 7.5 had got the
equivalence of (1) and $(1)^{op}$ for noetherian rings. However,
the proof here is rather different from that in [19]. Moreover, to
prove such an equivalence (Proposition 3.5), we get some other
results (for example, Lemmas 3.2 and 3.3), which are of
independent interest themselves. As corollaries of Theorem above,
we get a new characterization of $k$-Gorenstein algebras, and we
in addition have that for a faithfully balanced selforthogonal
bimodule $_{\Lambda}U _{\Lambda}$ its left injective dimension and
right injective dimension are identical provided $_{\Lambda}U$ (or
$U _{\Lambda}$) is Gorenstein. We in Section 4 study the dual
theory relative to Gorenstein modules (Theorems 4.1 and 4.4). In
the final section we give some properties of $_{\Lambda}U
_{\Gamma}$ with finite left or right injective dimension. Some
known results in [8] and [16] are obtained as corollaries.

\vspace{0.3cm}

\centerline{\bf 2. Preliminaries}

\vspace{0.1cm}

In this section we give some lemmas, which are useful in the rest
part of this paper.

Suppose that $A\in$mod $\Lambda$ (resp. mod $\Gamma ^{op}$). We
call Hom$_{\Lambda}(_{\Lambda}A, {_{\Lambda}U} _{\Gamma})$ (resp.
Hom$_{\Gamma}(A_{\Gamma}, {_{\Lambda}U} _{\Gamma})$) the dual
module of $A$ with respect to $U$, and denote either of these
modules by $A^*$. For a homomorphism $f$ between $\Lambda$-modules
(resp. $\Gamma ^{op}$-modules), we put $f^*={\rm Hom}(f,
{_{\Lambda}U} _{\Gamma})$. Let $\sigma _{A}: A \rightarrow A^{**}$
via $\sigma _{A}(x)(f)=f(x)$ for any $x \in A$ and $f \in A^*$ be
the canonical evaluation homomorphism. $A$ is called $U$-reflexive
if $\sigma _{A}$ is an isomorphism. Under the assumption of
$_{\Lambda}U _{\Gamma}$ being faithfully balanced, it is easy to
see that any projective module in mod $\Lambda$ (resp. mod
$\Lambda ^{op}$) is $U$-reflexive.

For a $\Lambda$-module (resp. $\Lambda ^{op}$-module) $X$, we use
{\it l.}fd$_{\Lambda}(X)$ (resp. {\it r.}fd$_{\Lambda}(X)$) to
denote the left (resp. right) flat dimension of $X$, and use {\it
l.}id$_{\Lambda}(X)$ (resp. {\it r.}id$_{\Lambda}(X)$) to denote
its left (resp. right) injective dimension. For a $\Lambda$-module
(resp. $\Gamma ^{op}$-module) $Y$, we denote either of
Hom$_{\Lambda}(_{\Lambda}U _{\Gamma}, {_{\Lambda}Y})$ and
Hom$_{\Gamma}(_{\Lambda}U _{\Gamma}, Y_{\Gamma})$ by $^*Y$.

\vspace{0.1cm}

{\bf Lemma 2.1} {\it Let} $\Lambda$ {\it and} $\Gamma$ {\it be
left and right noetherian rings and} $n$ {\it a non-negative
integer. If} $_{\Lambda}E$ ({\it resp.} $E_{\Gamma}$) {\it is
injective, then} {\it l.}fd$_{\Gamma}(^*E)$ ({\it resp.
r.}fd$_{\Lambda}(^*E))\leq n$ {\it if and only if}
Hom$_{\Lambda}($Ext$_{\Gamma}^{n+1}(A, U), E)$ ({\it resp.}
Hom$_{\Gamma}($Ext$_{\Lambda}^{n+1}(A, U), E))=0$ {\it for any} $A
\in$mod $\Gamma ^{op}$ ({\it resp.} mod $\Lambda$).

\vspace{0.2cm}

{\it Proof.} It is trivial by [6] Chapter VI, Proposition 5.3.
$\blacksquare$

\vspace{0.2cm}

We use add$_{\Lambda}U$ (resp. add$U_{\Gamma}$) to denote the
subcategory of mod $\Lambda$ (resp. mod $\Gamma ^{op}$) consisting
of all modules isomorphic to direct summands of finite sums of
copies of $_{\Lambda}U$ (resp. $U_{\Gamma}$). Let $A \in$mod
$\Lambda$. If there is an exact sequence $\cdots \to U _{n} \to
\cdots \to U _{1} \to U _{0} \to A \to 0$ in \linebreak mod
$\Lambda$ with each $U _{i}\in$add$_{\Lambda}U$ for any $i \geq
0$, then we define the $U$-resolution dimension of $A$, denoted by
$U$-resol.dim$_{\Lambda}(A)$, as inf$\{ n|$ there is an exact
sequence $0 \to U _{n} \to \cdots \to U _{1} \to U _{0} \to A \to
0$ in mod $\Lambda$ with each $U _{i}\in$add$_{\Lambda}U$ for any
$0 \leq i \leq n \}$. We set $U$-resol.dim$_{\Lambda}(A)$ infinity
if no such an integer exists. Similarly, for a module $B$ in mod
$\Gamma ^{op}$, we may define $U$-resol.dim$_{\Gamma}(B)$.

\vspace{0.2cm}

{\bf Lemma 2.2} {\it Let} $\Lambda$ {\it and} $\Gamma$ {\it be
left and right noetherian rings and} $n$ {\it a non-negative
integer. For a module} $X$ {\it in} mod $\Gamma ^{op}$, {\it if}
grade$_{U}X \geq n$ {\it and} grade$_{U}$Ext$_{\Gamma}^n(X, U)\geq
n+1$, {\it then} Ext$_{\Gamma}^{n}(X, U)=0$.

\vspace{0.2cm}

{\it Proof.} The proof of [13] Lemma 2.6 remains valid here, we
omit it. $\blacksquare$

\vspace{0.2cm}

{\bf Lemma 2.3} ([13] Lemma 2.7) {\it Let} $E \in$ mod $\Lambda$
({\it resp.} mod $\Gamma ^{op}$) {\it be injective. Then} {\it
l.}fd$_{\Gamma}(^*E)$ ({\it resp.} {\it r.}fd$_{\Lambda}(^*E))
\leq n$ {\it if and only if} $U$-resol.dim$_{\Lambda}(E)$ ({\it
resp.} $U$-resol.dim$_{\Gamma}(E)) \leq n$.

\vspace{0.2cm}

{\bf Lemma 2.4} ([13] Proposition 3.2) {\it The following
statements are equivalent.}

(1) $U$-dom.dim$(_{\Lambda}U) \geq 1$.

(2) $(-)^{**}:$ mod $\Lambda \to$ mod $\Lambda$ {\it preserves
monomorphisms}.

(3) $0 \to (_{\Lambda}U)^{**} \buildrel {f_{0}^{**}} \over
\longrightarrow E_{0}^{**}$ {\it is exact}.

(1)$^{op}$ $U$-dom.dim$(U _{\Gamma}) \geq 1$.

(2)$^{op}$ $(-)^{**}:$ mod $\Gamma ^{op} \to$ mod $\Gamma ^{op}$
{\it preserves monomorphisms}.

(3)$^{op}$ $0 \to (U _{\Gamma})^{**} \buildrel {(f'_{0})^{**}}
\over \longrightarrow (E'_{0})^{**}$ {\it is exact}.

\vspace{0.2cm}

{\bf Lemma 2.5} ([11] Lemma 2.7) {\it Let} $\Lambda$ {\it and}
$\Gamma$ {\it be left and right noetherian rings. The following
statements are equivalent.}

(1) $M^*$ {\it is} $U$-{\it reflexive for any} $M\in$ mod $\Lambda$.

(2) $[{\rm Ext}_{\Lambda}^2(M, U)]^*=0$ {\it for any} $M\in$ mod
$\Lambda$.

(1)$^{op}$ $N^*$ {\it is} $U$-{\it reflexive for any} $N\in$ mod
$\Gamma ^{op}$.

(2)$^{op}$ $[{\rm Ext}_{\Gamma}^2(N, U)]^*=0$ {\it for any} $N\in$
mod $\Gamma ^{op}$.

\vspace{0.2cm}

From now on, assume that

$$0 \to {_{\Lambda}U} \to E_{0} \to E_{1}
\cdots \to E_{i} \to \cdots$$ is a minimal injective resolution of
$_{\Lambda}U$, and

$$0 \to U _{\Gamma} \to E_{0}' \to E_{1}'
\cdots \to E_{i}' \to \cdots$$ is a minimal injective resolution
of $U _{\Gamma}$.

\vspace{0.2cm}

{\bf Lemma 2.6} ([14] Corollary 3.7) (1) $U$-resol.dim$_{\Lambda}
(E_{i})\leq i$ {\it for any} $0 \leq i \leq k-1$ {\it if and only
if} s.grade$_{U}{\rm Ext} _{\Gamma}^{i}(N, U) \geq i$ {\it for any}
$N\in$mod $\Gamma ^{op}$ {\it and} $1 \leq i \leq k$.

(2) $U$-resol.dim$_{\Gamma} (E_{i}')\leq i$ {\it for any} $0 \leq
i \leq k-1$ {\it if and only if} s.grade$_{U}{\rm Ext}
_{\Lambda}^{i}(M, U) \geq i$ {\it for any} $M\in$mod $\Lambda$
{\it and} $1 \leq i \leq k$.

\vspace{0.2cm}

{\bf Lemma 2.7} $_{\Lambda}U$ {\it is 1-Gorenstein if and only if}
$U _{\Gamma}$ {\it is 1-Gorenstein}.

\vspace{0.2cm}

{\it Proof.} By [13] Corollary 2.5, ${\rm
Hom}_{\Lambda}(_{\Lambda}U _{\Gamma}, E_0)$ is left $\Gamma$-flat
if and only if ${\rm Hom}_{\Gamma}(_{\Lambda}U _{\Gamma}, E'_0)$
is right $\Lambda$-flat. By Lemma 2.3, we then have that $E_0$ is
in add$_{\Lambda}U$ if and only if $E'_0$ is in add$U _{\Gamma}$.
So, it follows from Lemma 2.6 that s.grade$_{U}{\rm Ext}
_{\Gamma}^{1}(N, U) \geq 1$ for any $N\in$mod $\Gamma ^{op}$ if
and only if s.grade$_{U}{\rm Ext} _{\Lambda}^{1}(M, U) \geq 1$ for
any $M\in$mod $\Lambda$. Hence we conclude that $_{\Lambda}U$ is
1-Gorenstein if and only if $U _{\Gamma}$ is 1-Gorenstein.
$\blacksquare$

\vspace{0.5cm}

\centerline{\bf 3. Characterizations of $k$-Gorenstein modules}

\vspace{0.2cm}

In this section, we characterize $k$-Gorenstein modules in terms of
the $U$-resolution dimension of some special injective modules and
the property of the functors ${\rm Ext}^{i}({\rm Ext}^{i}(-, U), U)$
preserving monomorphisms, and also establish the left-right symmetry
of the notion of $k$-Gorenstein modules by using different methods
from that in [19]. In order to get our main theorem, we need some
lemmas.

\vspace{0.2cm}

{\bf Lemma 3.1} {\it If} $_{\Lambda}U$ {\it is} $k$-{\it
Gorenstein, then} Ext$_{\Lambda}^{i}($Ext$_{\Gamma}^{i}(-, U),
U):$ mod $\Gamma ^{op} \to$mod $\Gamma ^{op}$ {\it preserves
monomorphisms for any} $0 \leq i \leq k-1$.

\vspace{0.2cm}

{\it Proof.} We proceed by induction on $k$. The case for $k=1$
follows from Lemma 2.6 and Lemma 2.4.

Now suppose $k \geq 2$ and $0 \to X \to Y \to Z \to 0$ is an exact
sequence in mod $\Gamma ^{op}$. Then we have in mod $\Lambda$ the
following commutative diagram with the row exact:

$$\begin{tabular}{cccccccccc}
& ${\rm Ext}_{\Gamma}^{k-1}(Z, U)$ & $\buildrel {\alpha} \over
\longrightarrow$ & ${\rm Ext}_{\Gamma}^{k-1}(Y, U)$ & $\buildrel
{\beta} \over \longrightarrow$ & ${\rm Ext}_{\Gamma}^{k-1}(X, U)$
& $\buildrel {\gamma} \over \longrightarrow$
& ${\rm Ext}_{\Gamma}^{k}(Z, U),$ \\
& & $\searrow$ $ $ $\nearrow$ & & $\searrow$ $ $ $\nearrow$
& & $\searrow$ $ $ $\nearrow$ &\\
& & $A$ & & $B$ & & $C$ &
\end{tabular}$$
where $A={\rm Im}\alpha$, $B={\rm Im}\beta$ and $C={\rm
Im}\gamma$, and each triangle in above diagram is an epic-monic
resolution. Since $_{\Lambda}U$ is $k$-Gorenstein,
s.grade$_{U}{\rm Ext}^{i}_{\Gamma} (N, U)\geq i$ for any $N
\in$mod $\Gamma ^{op}$ and $1 \leq i \leq k$. So grade$_{U}A \geq
k-1$, grade$_{U}B \geq k-1$, grade$_{U}C \geq k$ and we have exact
sequences:

$$0={\rm Ext}_{\Lambda}^{k-1}(C, U) \longrightarrow
{\rm Ext}_{\Lambda}^{k-1}({\rm Ext}_{\Gamma}^{k-1}(X, U), U)
\longrightarrow {\rm Ext}_{\Lambda}^{k-1}(B, U),$$

$$0={\rm Ext}_{\Lambda}^{k-2}(A, U) \longrightarrow
{\rm Ext}_{\Lambda}^{k-1}(B, U)  \longrightarrow {\rm
Ext}_{\Lambda}^{k-1}({\rm Ext}_{\Gamma}^{k-1}(Y, U), U)$$ and we
then get a composition of monomorphisms:

$${\rm Ext}_{\Lambda}^{k-1}({\rm Ext}_{\Gamma}^{k-1}(X, U),
U) \hookrightarrow {\rm Ext}_{\Lambda}^{k-1}(B, U) \hookrightarrow
{\rm Ext}_{\Lambda}^{k-1}({\rm Ext}_{\Gamma}^{k-1}(Y, U), U),$$
which is also a monomorphism. $\blacksquare$

\vspace{0.2cm}

Let $M$ be in mod $\Lambda$ (resp. mod $\Gamma ^{op}$) and $P_{1}
\buildrel {f} \over \longrightarrow P_{0} \to M \to 0$ a projective
resolution of $M$ in mod $\Lambda$ (resp. mod $\Gamma ^{op}$). Then
we have an exact sequence $0 \to M^* \to P_{0}^* \buildrel {f^*}
\over \longrightarrow P_{1}^* \to X \to 0$, where $X={\rm
Coker}f^*$. For a positive integer $k$, recall from [10] that $M$ is
called $U$-$k$-torsionfree if Ext$_{\Gamma}^{i} (X, U)=0$ (resp.
Ext$_{\Lambda}^{i} (X, U)=0)$ for any $1 \leq i \leq k$. $M$ is
called $U$-$k$-syzygy if there is an exact sequence $0 \to M \to
X_{0} \to X_{1} \to \cdots \to X_{k-1}$ with all $X_{i}$ in
add$_{\Lambda}U$ (resp. add$U_{\Gamma}$). Put $_{\Lambda}U
_{\Gamma}={_{\Lambda}\Lambda} _{\Lambda}$, then, in this case, the
notions of $U$-$k$-torsionfree modules and $U$-$k$-syzygy modules
are just that of $k$-torsionfree modules and $k$-syzygy modules
respectively (see [2] for the definitions of $k$-torsionfree modules
and $k$-syzygy modules). We use $\cal{T}$$_{U}^{k}({\rm mod}\
\Lambda)$ (resp. $\cal{T}$$_{U}^{k}({\rm mod}\ \Gamma ^{op})$) and
$\Omega ^{k}_{U}({\rm mod}\ \Lambda)$ (resp. $\Omega ^{k}_{U}({\rm
mod}\ \Gamma ^{op})$) to denote the full subcategory of mod
$\Lambda$ (resp. mod $\Gamma ^{op}$) consisting of
$U$-$k$-torsionfree modules and $U$-$k$-syzygy modules,
respectively. It was in [10] pointed out that
$\cal{T}$$_{U}^{k}({\rm mod}\ \Lambda) \subseteq \Omega
^{k}_{U}({\rm mod}\ \Lambda)$ and $\cal{T}$$_{U}^{k}({\rm mod}\
\Gamma ^{op}) \subseteq \Omega ^{k}_{U}({\rm mod}\ \Gamma ^{op})$.

The following two lemmas are of independent interest themselves.

\vspace{0.2cm}

{\bf Lemma 3.2} {\it Let} $\Lambda$ {\it and} $\Gamma$ {\it be left
and right noetherian rings. If} grade$_{U}{\rm Ext}^{i+1}_{\Lambda}
(M, U)\geq i$ {\it for any} $M \in$mod $\Lambda$ {\it and} $1 \leq i
\leq k-1$, {\it then each} $k$-{\it syzygy module in} mod $\Lambda$
{\it is in} $\Omega ^{k}_{U}({\rm mod}\ \Lambda)$.

\vspace{0.2cm}

{\it Proof.} Suppose that grade$_{U}{\rm Ext}^{i+1}_{\Lambda} (M,
U)\geq i$ for any $M \in$mod $\Lambda$ and $1 \leq i \leq k-1$. Then
by [11] Theorem 3.1, $\Omega ^{k}_{U}({\rm mod}\
\Lambda)=\cal{T}$$_{U}^{k}({\rm mod}\ \Lambda)$. So it suffices to
show that each $k$-syzygy module in mod $\Lambda$ is in
$\cal{T}$$_{U}^{k}({\rm mod}\ \Lambda)$. The following proof is
similar to that of $(1)\Rightarrow (2)$ in [11] Theorem 3.1. For the
sake of completeness, we give here the proof.

We proceed by induction on $k$.

Notice that $\Lambda$ is $U$-reflexive, it follows easily that each
1-syzygy module in mod $\Lambda$ is in $\Omega ^{1}_{U}({\rm mod}\
\Lambda)(=\cal{T}$$_{U}^{1}({\rm mod}\ \Lambda))$.

Assume that $k=2$ and $M$ is a 2-syzygy module in mod $\Lambda$.
Then there is an exact sequence $0 \to M \to P_1 \buildrel {f} \over
\to P_0$ in mod $\Lambda$ with $P_0$ and $P_1$ projective. By [11]
Lemma 2.4, $M \cong ({\rm Coker}f^*)^*$. It follows from Lemma 2.5
and [10] Lemma 4 that $M$ is $U$-reflexive and $U$-2-torsionfree.
The case for $k=2$ follows.

Now suppose that $k \geq 3$ and $M$ is a $k$-syzygy module in mod
$\Lambda$. Then there is an exact sequence: $$P_{k+1} \buildrel
{f_{k+1}} \over \longrightarrow P_k \buildrel {f_k} \over
\longrightarrow P_{k-1} \buildrel {f_{k-1}} \over \longrightarrow
\cdots \buildrel {f_2} \over \longrightarrow P_1 \buildrel {f_1}
\over \longrightarrow P_0 \to X \to 0$$ in mod $\Lambda$ such that
$M={\rm Coker}f_{k+1}$, where each $P_i$ is projective for any $0
\leq i \leq k+1$. By induction assumption, $M \in
\cal{T}$$_{U}^{k-1}({\rm mod}\ \Lambda)$. We will show that $M \in
\cal{T}$$_{U}^{k}({\rm mod}\ \Lambda)$. Notice that $k \geq 3$, so
$M$ is $U$-reflexive and hence it suffices to show that
Ext$_{\Gamma}^i(M^*, U)=0$ for any $1 \leq i \leq k-2$ by [11] Lemma
2.9.

Put $N=$Coker$f_{k-1}^{*}$. Then, by [11] Lemma 2.4, $M \cong N^{*}$
and $M^{*} \cong N^{**}$. We claim that Ext$_{\Gamma}^{i}(N, U)=0$
for any $1 \leq i \leq k-2$. If $k=3$, then Coker$f_{k-1}$ is a
submodule of $P_{0}$. But $P_{0}$ is $U$-reflexive, so
Coker$f_{k-1}$ is $U$-torsionless. By [15] Lemma 2.1,
Ext$_{\Gamma}^{1}(N, U) \cong$Ker$\sigma _{{\rm Coker}f_{k-1}}=0$.
If $k=4$, then Coker$f_{k-1}$ is a 2-syzygy module in mod $\Lambda$
and so Coker$f_{k-1}$ is $U$-reflexive by the above argument. Thus
by [15] Lemma 2.1, Ext$_{\Gamma}^{1}(N, U) \cong$Ker$\sigma _{{\rm
Coker}f_{k-1}}=0$ and Ext$_{\Gamma}^{2}(N, U) \cong$Coker$\sigma
_{{\rm Coker}f_{k-1}}=0$ and the case for $k=4$ follows. If $k \geq
5$, then Coker$f_{k-1}$ is a $(k-2)$-syzygy module in mod $\Lambda$
and so Coker$f_{k-1} \in \cal{T}$$_{U}^{k-2}({\rm mod}\ \Lambda)$ by
induction assumption. It is clear that ${\rm Coker}f_{k-1}$ is
$U$-reflexive. Then by using an argument similar to that in the
proof of the case for $k=4$, we have that Ext$_{\Gamma}^{1}(N,
U)=0=$Ext$_{\Gamma}^{2}(N, U)=0$. On the other hand, by [11] Lemma
2.9, we have that Ext$_{\Gamma}^{i}(({\rm Coker}f_{k-1})^{*}, U)=0$
for any $1 \leq i \leq k-4$. It follows from the exact sequence $0
\to ({\rm Coker}f_{k-1})^{*} \to P_{k-2}^{*} \buildrel {f_{k-1}^{*}}
\over \longrightarrow P_{k-1}^{*} \to N \to 0$ that
Ext$_{\Gamma}^{i}(N, U)=0$ for any $3 \leq i \leq k-2$. So
Ext$_{\Gamma}^{i}(N, U)=0$ for any $1 \leq i \leq k-2$.

By [15] Lemma 2.1, we have an exact sequence: $$0 \to {\rm
Ext}_{\Lambda}^{1}({\rm Coker}f_{k-1}, U) \to N \buildrel {\sigma
_{N}} \over \longrightarrow N^{**} \to {\rm Ext}_{\Lambda}^{2}({\rm
Coker}f_{k-1}, U) \to 0.$$ Then Ker$\sigma _N \cong {\rm
Ext}_{\Lambda}^{1}({\rm Coker}f_{k-1}, U) \cong {\rm
Ext}_{\Lambda}^{k-1}(X, U)$ and Coker$\sigma _N \cong {\rm
Ext}_{\Lambda}^{2}({\rm Coker}f_{k-1}, U) \cong {\rm
Ext}_{\Lambda}^{k}(X, U)$. So we get the following exact sequences:
$$0 \to {\rm Ext}_{\Lambda}^{k-1}(X, U) \to N
\buildrel {\pi} \over \longrightarrow {\rm Im}\sigma _{N} \to 0
\eqno{(1)}$$
$$0 \to {\rm Im}\sigma _{N} \buildrel {\mu} \over \longrightarrow
N^{**} \to {\rm Ext}_{\Lambda}^{k}(X, U) \to 0 \eqno{(2)}$$ where
$\sigma _{N}=\mu \pi$. Since Ext$_{\Gamma}^{i}(N, U)=0$ for any $1
\leq i \leq k-2$ and grade$_{U}$Ext$_{\Lambda}^{k-1}(X, U) \geq
k-2$, from the exact sequence (1) we have Ext$_{\Gamma}^{i}({\rm
Im}\sigma _N, U)=0$ for any $1 \leq i \leq k-2$. Moreover, since
grade$_{U}$Ext$_{\Lambda}^{k}(X, U) \geq k-1$, from the exact
sequence (2) we get that Ext$_{\Gamma}^{i}(N^{**}, U)=0$ for any $1
\leq i \leq k-2$, which yields Ext$_{\Gamma}^{i}(M^{*}, U)=0$ for
any $1 \leq i \leq k-2$. We are done. $\blacksquare$

\vspace{0.2cm}

{\bf Lemma 3.3} {\it Let} $\Lambda$ {\it and} $\Gamma$ {\it be
left and right noetherian rings. For a positive integer} $k$, {\it
the following statements are equivalent}.

(1) grade$_{U}$Ext$_{\Lambda}^{i+1}(M, U)\geq i$ {\it for any}
$M\in$ mod $\Lambda$ {\it and} $1\leq i \leq k-1$.

(2) Ext$_{\Gamma}^{i-1}($Ext$_{\Lambda}^{i+1}(M, U), U)=0$ {\it
for any} $M\in$ mod $\Lambda$ {\it and} $1\leq i \leq k-1$.

(1)$^{op}$ grade$_{U}$Ext$_{\Gamma}^{i+1}(N, U)\geq i$ {\it for
any} $N\in$ mod $\Gamma ^{op}$ {\it and} $1\leq i \leq k-1$.

(2)$^{op}$ Ext$_{\Lambda}^{i-1}($Ext$_{\Gamma}^{i+1}(N, U), U)=0$
{\it for any} $N\in$ mod $\Gamma ^{op}$ {\it and} $1\leq i \leq
k-1$.

\vspace{0.2cm}

{\it Proof.} The implications that $(1)\Rightarrow (2)$ and
$(1)^{op}\Rightarrow (2)^{op}$ are trivial.

$(2)\Rightarrow (1)$ We proceed by induction on $k$. It is trivial
when $k=1$ or $k=2$. Now suppose $k\geq 3$. By induction
assumption, for any $M\in$ mod $\Lambda$, we have
grade$_{U}$Ext$_{\Lambda}^{i+1}(M, U)\geq i$ for any $1\leq i \leq
k-2$ and grade$_{U}$Ext$_{\Lambda}^{k}(M, U)\geq k-2$. In
addition, by (2), we have
Ext$_{\Gamma}^{k-2}($Ext$_{\Lambda}^{k}(M, U), U)=0$. So
grade$_{U}$Ext$_{\Lambda}^{k}(M, U)\geq k-1$.

$(1)^{op} \Rightarrow (1)$ We also proceed by induction on $k$.
The case $k=1$ is trivial. The case $k=2$ follows from Lemma 2.5.
Now suppose $k\geq 3$.

Let $M\in$ mod $\Lambda$ and
$$\cdots \to P_i \to \cdots \to P_1 \to P_0 \to M \to 0$$ a
projective resolution of $M$ in mod $\Lambda$. Put
$M_i=$Coker$(P_i \to P_{i-1})$ (where $M_1=M$) and
$X_i=$Coker$(P_{i-1}^* \to P_i^*)$ for any $i\geq 1$. By induction
assumption, we have grade$_{U}$Ext$_{\Lambda}^{i+1}(M, U)$
\linebreak $\geq i$ for any $1\leq i \leq k-2$ and
grade$_{U}$Ext$_{\Lambda}^k(M, U)\geq k-2$. So it suffices to
prove \linebreak Ext$_{\Gamma}^{k-2}$(Ext$_{\Lambda}^k(M, U),
U)=0$.

By [11] Theorem 3.1, $\Omega _{U} ^i({\rm mod}\
\Lambda)=\mathcal{T} _{U}^i({\rm mod}\ \Lambda)$ for any $1\leq i
\leq k-1$. For any $t\geq k$, since $M_t \in \Omega _{U}
^{k-1}({\rm mod}\ \Lambda)$ by Lemma 3.2, $M_t \in \mathcal{T}
_{U}^{k-1}({\rm mod}\ \Lambda)$. It follows that
Ext$_{\Gamma}^{j}(X_t, U)=0$ for any $1 \leq j \leq k-1$ and $t
\geq k$.

On the other hand, by [12] Lemma 2 we have an exact sequence:
$$0\to {\rm Ext}_{\Lambda}^k(M, U) \to X_k \to P_{k+1}^* \to
X_{k+1} \to 0.$$ Put $K=$Im$(X_k \to P_{k+1}^*)$. From the
exactness of $0\to K \to P_{k+1}^* \to X_{k+1} \to 0$ we know that
Ext$_{\Gamma}^j(K, U)=0$ for any $1 \leq j \leq k-2$ and ${\rm
Ext}_{\Gamma}^k(X_{k+1}, U)\cong {\rm Ext}_{\Gamma}^{k-1}(K, U)$.
Moreover, from the exactness of $0\to {\rm Ext}_{\Lambda}^k(M, U)
\to X_k \to K \to 0$ we know that ${\rm Ext}_{\Gamma}^{k-1}(K,
U)\cong {\rm Ext}_{\Gamma}^{k-2}({\rm Ext}_{\Lambda}^k(M, U), U)$.
So ${\rm Ext}_{\Gamma}^{k-2}({\rm Ext}_{\Lambda}^k(M, U), U)\cong
{\rm Ext}_{\Gamma}^k(X_{k+1}, U)$. By (1)$^{op}$, we then have
that grade$_{U}{\rm Ext}_{\Gamma}^{k-2}({\rm Ext}_{\Lambda}^k(M,
U), U)=$grade$_{U}{\rm Ext}_{\Gamma}^k(X_{k+1}, U)\geq k-1$. It
follows from Lemma 2.2 that
Ext$_{\Lambda}^{k-2}$(Ext$_{\Lambda}^k(M, U), U)=0$.

By symmetry, we have the implications of $(2)^{op} \Rightarrow
(1)^{op}$ and $(1) \Rightarrow (1)^{op}$. We are done.
$\blacksquare$

\vspace{0.2cm}

The following result not only generalizes [2] Proposition 2.26,
but also means that the statements in this proposition are
left-right symmetric.

\vspace{0.2cm}

{\bf Corollary 3.4} {\it Let} $\Lambda$ {\it be a left and right
noetherian ring. For a positive integer} $k$, {\it the following
statements are equivalent}.

(1) gradeExt$_{\Lambda}^{i+1}(M, \Lambda)\geq i$ {\it for any}
$M\in$ mod $\Lambda$ {\it and} $1\leq i \leq k-1$.

(2) Ext$_{\Lambda}^{i-1}($Ext$_{\Lambda}^{i+1}(M, \Lambda),
\Lambda)=0$ {\it for any} $M\in$ mod $\Lambda$ {\it and} $1\leq i
\leq k-1$.

(3) $\Omega _{\Lambda}^i({\rm mod}\ \Lambda)=\mathcal{T}
_{\Lambda}^i({\rm mod}\ \Lambda)$ {\it for any} $1\leq i \leq k$.

(1)$^{op}$ gradeExt$_{\Lambda}^{i+1}(N, \Lambda)\geq i$ {\it for
any} $N\in$ mod $\Lambda ^{op}$ {\it and} $1\leq i \leq k-1$.

(2)$^{op}$ Ext$_{\Lambda}^{i-1}($Ext$_{\Lambda}^{i+1}(N, \Lambda),
\Lambda)=0$ {\it for any} $N\in$ mod $\Lambda ^{op}$ {\it and}
$1\leq i \leq k-1$.

(3)$^{op}$ $\Omega  _{\Lambda}^i({\rm mod}\ \Lambda
^{op})=\mathcal{T} _{\Lambda}^i({\rm mod}\ \Lambda ^{op})$ {\it
for any} $1\leq i \leq k$.

\vspace{0.2cm}

{\it Proof.} By Lemma 3.3 we have $(1) \Leftrightarrow (2)
\Leftrightarrow (1)^{op} \Leftrightarrow (2)^{op}$, and by [2]
Proposition 2.26 we have $(1) \Leftrightarrow (3)$ and $(1)^{op}
\Leftrightarrow (3)^{op}$. $\blacksquare$

\vspace{0.2cm}

The following proposition, had been got by Wakamatsu in [19]
Theorem 7.5 for noetherian rings, shows the left-right symmetry of
the notion of $k$-Gorenstein modules. However, the proof here is
rather different from that in [19].

\vspace{0.2cm}

{\bf Proposition 3.5} {\it For a positive integer} $k$, {\it the
following statements are equivalent}.

(1) s.grade$_{U}{\rm Ext}^{i}_{\Lambda} (M, U)\geq i$ {\it for
any} $M \in$mod $\Lambda$ {\it and} $1 \leq i \leq k$.

(2) s.grade$_{U}{\rm Ext}^{i}_{\Gamma} (N, U)\geq i$ {\it for any}
$N \in$mod $\Gamma ^{op}$ {\it and} $1 \leq i \leq k$.

\vspace{0.2cm}

{\it Proof.} By symmetry, we only need to prove (2) implies (1).

We proceed by induction on $k$. The case $k=1$ follows from Lemma
2.7. Now suppose that $k\geq 2$ and s.grade$_{U}{\rm
Ext}^{i}_{\Gamma} (N, U)\geq i$ for any $N \in$mod $\Gamma ^{op}$
and $1 \leq i \leq k$. By Lemma 3.3 and [11] Theorem 3.1, we have
that $\cal{T}$$^{i}_{U}({\rm mod}\ \Lambda)=\Omega ^{i}_{U}({\rm
mod}\ \Lambda)$ for any $1 \leq i \leq k$. By induction
assumption, for any $M \in$mod $\Lambda$ we have that
s.grade$_{U}{\rm Ext}^{i}_{\Lambda} (M, U)\geq i$ for any $1 \leq
i \leq k-1$ and s.grade$_{U}{\rm Ext}^{k}_{\Lambda} (M, U)\geq
k-1$.

Assume that
$$\cdots \buildrel {f_{k}} \over \longrightarrow
P_{k-1} \buildrel {f_{k-1}} \over \longrightarrow \cdots
\longrightarrow P_{1} \buildrel {f_{1}} \over \longrightarrow
P_{0} \longrightarrow M \longrightarrow 0$$ is a (minimal)
projective resolution of $M$ in mod $\Lambda$. For any $t \geq k$,
since Coker$f_{t}$ is $(k-1)$-syzygy, by Lemma 3.2 we have that
Coker$f_{t}\in \Omega ^{k-1}_{U}({\rm mod}\ \Lambda)$ and
Coker$f_{t}\in \cal{T}$$^{k-1}_{U}({\rm mod}\ \Lambda)$, that is,
Coker$f_{t}$ is $U$-$(k-1)$-torsionfree, which implies that
Ext$_{\Gamma}^{j}({\rm Coker}f^*_{t}, U)=0$ for any $1 \leq j \leq
k-1$ and $t \geq k$.

Let $X$ be a submodule of ${\rm Ext}_{\Lambda}^{k}(M, U)$. Then
grade$_{U}X \geq k-1$. On the other hand, by [15] Lemma 2.1 there
is an exact sequence $0 \to {\rm Ext}_{\Lambda}^{k}(M, U)\to {\rm
Coker}f^*_{k} \buildrel {\sigma _{{\rm Coker}f^*_k}} \over
\longrightarrow ({\rm Coker}f^*_{k})^{**} \to {\rm
Ext}_{\Lambda}^{k+1}(M, U)\to 0$ and then there is a composition
of monomorphisms: $X \hookrightarrow {\rm Ext}_{\Lambda}^{k}(M, U)
\hookrightarrow {\rm Coker}f_{k}^*$. Put $Y=$Coker$(X
\hookrightarrow {\rm Coker}f_k^*)$. Notice that
Ext$_{\Gamma}^{j}({\rm Coker}f^*_{k}, U)=0$ for any $1 \leq j \leq
k-1$, we then have an embedding Ext$_{\Gamma}^{k-1}(X,
U)\hookrightarrow$Ext$_{\Gamma}^k(Y, U)$. By assumption,
s.grade$_{U}$Ext$_{\Gamma}^k(Y, U)\geq k$. So grade$_{U}{\rm
Ext}_{\Gamma}^{k-1}(X, U) \geq k$ and hence ${\rm
Ext}_{\Gamma}^{k-1}(X, U)=0$ by Lemma 2.2. It follows that
grade$_{U}X \geq k$ and s.grade$_{U}{\rm Ext}_{\Lambda}^k(M, U)
\geq k$. We are done. $\blacksquare$

\vspace{0.2cm}

{\bf Lemma 3.6} {\it If} Ext$_{\Gamma}^{i}($Ext$_{\Lambda}^{i}(-,
U), U)$: mod $\Lambda\to$mod $\Lambda$ {\it preserves
monomorphisms for any} $0 \leq i \leq k-1$, {\it then}
$_{\Lambda}U$ {\it is} $k$-{\it Gorenstein}.

\vspace{0.2cm}

{\it Proof.} We proceed by induction on $k$.

Assume that $(-)^{**}:$ mod $\Lambda \to$ mod $\Lambda$ preserves
monomorphisms, then, by Lemma 2.4, $U$-dom.dim$(_{\Lambda}U) \geq
1$ and $E_0$ is cogenerated by $_{\Lambda}U$. But $E_0$ is
finitely cogenerated, so $E_0 \in$add$_{\Lambda}U$. By Lemma 2.6,
we then have that s.grade$_{U}{\rm Ext} _{\Gamma}^{1}(N, U) \geq
1$ for any $N\in$mod $\Gamma ^{op}$ and $_{\Lambda}U$ is
1-Gorenstein. The case for $k=1$ is proved.

Now suppose $k\geq 2$ and
$$\cdots \buildrel {g_{k}} \over \longrightarrow
Q_{k-1} \buildrel {g_{k-1}} \over \longrightarrow \cdots
\longrightarrow Q_{1} \buildrel {g_{1}} \over \longrightarrow
Q_{0} \longrightarrow N \longrightarrow 0$$ is a (minimal)
projective resolution of a module $N$ in mod $\Gamma ^{op}$. By
induction hypothesis, $_{\Lambda}U$ is $(k-1)$-Gorenstein and
s.grade$_{U}{\rm Ext}_{\Gamma}^{i}(B, U) \geq i$ for any $B
\in$mod $\Gamma ^{op}$ and $1 \leq i \leq k-1$ (and certainly,
s.grade$_{U}{\rm Ext}_{\Gamma}^{i+1}(B, U) \geq i$ for any $B
\in$mod $\Gamma ^{op}$ and $1 \leq i \leq k-1$). So
$\cal{T}$$^{i}_{U}({\rm mod}\ \Gamma ^{op})=\Omega ^{i}_{U}({\rm
mod}\ \Gamma ^{op})$ for any $1 \leq i \leq k$ by [11] Theorem
3.1$^{op}$. By using a similar argument of $(2) \Rightarrow (1)$
in Proposition 3.5, we have that Ext$_{\Lambda}^{j}({\rm
Coker}g^*_{k}, U)=0$ for any $1 \leq j \leq k-1$.

Let $X$ be a submodule of ${\rm Ext}_{\Gamma}^{k}(N, U)$. Then
grade$_{U}X \geq k-1$. By using a similar argument of $(2)
\Rightarrow (1)$ in Proposition 3.5, we have a monomorphism $X
\hookrightarrow {\rm Coker}g_{k}^*.$ By assumption, $0 \to {\rm
Ext}_{\Gamma}^{k-1}({\rm Ext}_{\Lambda}^{k-1}(X, U), U) \to {\rm
Ext}_{\Gamma}^{k-1}({\rm Ext}_{\Lambda}^{k-1}({\rm Coker}g^*_{k},
U), U)(=0)$ is exact, so ${\rm Ext}_{\Gamma}^{k-1}({\rm
Ext}_{\Lambda}^{k-1}(X, U), U)=0$. On the other hand, by
Proposition 3.5, we have that s.grade$_{U}$Ext$_{\Lambda}^{k-1}(X,
U)\geq k-1$. So we conclude that grade$_{U}{\rm
Ext}_{\Lambda}^{k-1}(X, U) \geq k$ and hence ${\rm
Ext}_{\Lambda}^{k-1}(X, U)=0$ by Lemma 2.2. It follows that
grade$_{U}X \geq k$ and s.grade$_{U}{\rm Ext}_{\Gamma}^k(N, U)
\geq k$. We are done. $\blacksquare$

\vspace{0.2cm}

We are now in a position to state the main result in this paper.

\vspace{0.2cm}

{\bf Theorem 3.7} {\it The following statements are equivalent.}

(1) $_{\Lambda}U$ {\it is} $k$-{\it Gorenstein.}

(2) $U$-resol.dim$_{\Lambda}(E_{i}) \leq i$ {\it for any} $0 \leq
i \leq k-1$.

(3) Ext$_{\Gamma}^{i}({\rm Ext}_{\Lambda}^{i}(-, U), U)$: mod
$\Lambda \to$ mod $\Lambda$ {\it preserves monomorphisms} {\it for
any} $0 \leq i \leq k-1$.

(1)$^{op}$ $U _{\Gamma}$ {\it is} $k$-{\it Gorenstein.}

(2)$^{op}$ $U$-resol.dim$_{\Gamma}(E^{'}_{i}) \leq i$ {\it for
any} $0 \leq i \leq k-1$.

(3)$^{op}$ Ext$_{\Lambda}^{i}({\rm Ext}_{\Gamma}^{i}(-, U), U)$:
mod $\Gamma ^{op} \to$ mod $\Gamma ^{op}$ {\it preserves
monomorphisms} {\it for any} $0 \leq i \leq k-1$.

\vspace{0.2cm}

{\it Proof.} $(2)\Leftrightarrow (1) \Leftrightarrow (1)^{op}$ See
Lemma 2.6 and Proposition 3.5.

$(1)\Rightarrow (3)^{op}$ By Lemma 3.1.

$(3)\Rightarrow (1)$ By Lemma 3.6.

Symmetrically we have that $(2)^{op}\Leftrightarrow (1)^{op}$,
$(1)^{op}\Rightarrow (3)$ and $(3)^{op}\Rightarrow (1)^{op}$. The
proof is finished. $\blacksquare$

\vspace{0.2cm}

Put $_{\Lambda}U _{\Gamma}={_{\Lambda}\Lambda _{\Lambda}}$, by
Theorem 3.7, we then immediately have the following corollary, which
extends [7] Auslander's Theorem 3.7.

\vspace{0.2cm}

{\bf Corollary 3.8} {\it The following statements are equivalent.}

(1) s.grade$_{\Lambda}{\rm Ext}^{i}_{\Lambda}(M, \Lambda)\geq i$
{\it for any} $M \in$mod $\Lambda$ {\it and} $1 \leq i \leq k$.

(2) {\it The left flat dimension of the} $i$th {\it term of a
minimal injective resolution of} $_{\Lambda}\Lambda$ {\it is at
most} $i-1$ {\it for any} $1 \leq i \leq k$.

(3) Ext$_{\Lambda}^{i}({\rm Ext}_{\Lambda}^{i}(-, \Lambda),
\Lambda)$: mod $\Lambda \to$ mod $\Lambda$ {\it preserves
monomorphisms} {\it for any} $0 \leq i \leq k-1$.

(1)$^{op}$ s.grade$_{\Lambda}{\rm Ext}^{i}_{\Lambda}(N,
\Lambda)\geq i$ {\it for any} $N \in$mod $\Lambda ^{op}$ {\it and}
$1 \leq i \leq k$.

(2)$^{op}$ {\it The right flat dimension of the} $i$th {\it term
in a minimal injective resolution of} $\Lambda _{\Lambda}$ {\it at
most} $i-1$ {\it for any} $1 \leq i \leq k$.

(3)$^{op}$ Ext$_{\Lambda}^{i}({\rm Ext}_{\Lambda}^{i}(-, \Lambda),
\Lambda)$: mod $\Lambda ^{op} \to$ mod $\Lambda ^{op}$ {\it
preserves monomorphisms} {\it for any} $0 \leq i \leq k-1$.

\vspace{0.2cm}

$_{\Lambda}U _{\Gamma}$ is called a cotilting bimodule if
$_{\Lambda}U$ and $U _{\Gamma}$ are cotilting, that is, {\it
l.}id$_{\Lambda}(U)$ and {\it r.}id$_{\Gamma}(U)$ are
finite$^{[9]}$. If $_{\Lambda}U _{\Gamma}$ is a cotilting
bimodule, then {\it l.}id$_{\Lambda}(U)=${\it r.}id$_{\Gamma}(U)$
(see [4] Lemma 1.7). However, in general, we don't know whether
{\it l.}id$_{\Lambda}(U)<\infty$ implies that {\it
r.}id$_{\Gamma}(U)<\infty$. In fact, Auslander and Reiten in
[3]p.150 posed an important question which remains open: for an
artin algebra $\Lambda$, does {\it
l.}id$_{\Lambda}(\Lambda)<\infty$ imply {\it
r.}id$_{\Lambda}(\Lambda)<\infty$? Put
$_{\Lambda}U_{\Gamma}={_{\Lambda}U_{\Lambda}}$, as applications to
the results obtained above we have the following corollaries.

\vspace{0.2cm}

{\bf Corollary 3.9} {\it For a positive integer} $k$, {\it if
r.}id$_{\Lambda}(U)=k$ {\it and} $U_{\Lambda}$ {\it is}
$(k-1)$-{\it Gorenstein, then l.}id$_{\Lambda}(U)=k$.

\vspace{0.2cm}

{\it Proof.} The case for $k=1$ follows from [12] Corollary 1. Now
assume that {\it r.}id$_{\Lambda}(U)=k(\geq 2)$ and $U_{\Lambda}$
is $(k-1)$-Gorenstein. Then, by Lemma 2.6, we have
s.grade$_{U}{\rm Ext} _{\Lambda}^{k-1}(M, U) \geq k-1$ for any
$M\in$mod $\Lambda$. It follows from [12] Theorem that {\it
l.}id$_{\Lambda}(U) \leq 2k-2$. So {\it l.}id$_{\Lambda}(U)=k$ by
[4] Lemma 1.7. $\blacksquare$

\vspace{0.2cm}

{\bf Corollary 3.10} {\it l.}id$_{\Lambda}(U)$={\it
r.}id$_{\Lambda}(U)$ {\it if} $_{\Lambda}U$ ({\it or} $U
_{\Lambda}$) {\it is Gorenstein.}

\vspace{0.2cm}

{\it Proof.} By Theorem 3.7, Corollary 3.9 and its dual result.
$\blacksquare$

\vspace{0.2cm}

{\bf Corollary 3.11} ([5] Corollary 5.5) {\it Let} $\Lambda$ {\it
be a} $k$-{\it Gorenstein algebra for all} $k$. {\it Then} {\it
l.}id$_{\Lambda}(\Lambda)=$ {\it r.}id$_{\Lambda}(\Lambda)$.

\vspace{0.5cm}

\centerline{\bf 4. Dual theory}

\vspace{0.2cm}

In this section we study the dual theory relative to Gorenstein
bimodules.

For a non-negative integer $g$, we use $\mathcal{G}_{g}({\rm mod}\
\Lambda)$ (resp. $\mathcal{G}_{g}({\rm mod}\ \Gamma ^{op})$) to
denote the subcategory of mod $\Lambda$ (resp. mod $\Gamma ^{op}$)
consisting of the modules $M$ with grade$_{U}M=g$, and
$\mathcal{H}_{g}({\rm mod}\ \Lambda)$ (resp. $\mathcal{H}_{g}({\rm
mod}\ \Gamma ^{op})$) to denote the subcategory of
$\mathcal{G}_{g}({\rm mod}\ \Lambda)$ (resp. $\mathcal{G}_{g}({\rm
mod}\ \Gamma ^{op})$) consisting of the modules $M$ with ${\rm
Ext}_{\Lambda}^{i}(M, U)=0$ (resp. ${\rm Ext} _{\Gamma}^{i}(M,
U)=0$) for any $i\neq {\rm grade}_{U}M(=g)$.

\vspace{0.2cm}

{\bf Theorem 4.1} {\it Let} $\Lambda$ {\it and} $\Gamma$ {\it be
left and right noetherian rings and} $_{\Lambda}U _{\Gamma}$ {\it
a Gorenstein bimodule}.

(1) {\it If r.}id$_{\Gamma}(U)=g$, {\it then, for any} $0\neq M
\in \mathcal{G}_{g}({\rm mod}\ \Lambda)$, $M\cong {\rm
Ext}_{\Gamma}^{g}({\rm Ext}_{\Lambda}^{g}(M, U), U)$.

(2) {\it If r.}id$_{\Gamma}(U)=${\it l.}id$_{\Lambda}(U)=g$, {\it
then there is a duality between} $\mathcal{G}_{g}({\rm mod}\
\Lambda)$ {\it and} $\mathcal{G}_{g}({\rm mod}\ \Gamma ^{op})$
{\it given by} $M \to {\rm Ext}_{\Lambda}^{g}(M, U)$.

\vspace{0.2cm}

{\it Proof.} (1) Let $M$ be a non-zero module in
$\mathcal{G}_{g}({\rm mod}\ \Lambda)$ and
$$\cdots \to P_{i} \buildrel {f_{i}} \over \longrightarrow
\cdots \buildrel {f_{2}} \over \longrightarrow P_{1} \buildrel
{f_{1}} \over \longrightarrow P_{0} \to M \to 0$$ a projective
resolution of $M$ in mod $\Lambda$. If $g=0$, it is easy to see from
[15] Lemma 2.1 that $M \cong M^{**}$. Now suppose $g \geq 1$. Then
we get an exact sequence:

$$0 \to P_{0}^*
\buildrel {f_{1}^*} \over \longrightarrow P_{1}^* \buildrel
{f_{2}^*} \over \longrightarrow \cdots \buildrel {f_{g-1}^*} \over
\longrightarrow P_{g-1}^* \to ({\rm Im}f_{g})^* \to {\rm
Ext}_{\Lambda}^{g}(M, U) \to 0.$$

Since $_{\Lambda}U _{\Gamma}$ is a Gorenstein bimodule,
s.grade$_{U}{\rm Ext}_{\Lambda}^{g}(M, U)\geq g$. For any $i\geq
1$, put $K_{i}={\rm Coker}f_{i}^*$. We then have an exact
sequence:
$$P_{1}^{**}
\buildrel {f_{1}^{**}} \over \longrightarrow P_{0}^{**} \to {\rm
Ext}_{\Gamma}^{1}(K_{1}, U)(\cong {\rm
Ext}_{\Gamma}^{g-1}(K_{g-1}, U))\to 0.$$ On the other hand, we
have an exact sequence $0\to ({\rm Im}f_{g})^* \to P_{g}^*
\buildrel {f_{g+1}^*} \over \longrightarrow P_{g+1}^* \to K_{g+1}
\to 0$. Since {\it r.}id$_{\Gamma}(U)=g$, for any $i\geq g-1$ we
have ${\rm Ext}_{\Gamma}^{i}(({\rm Im}f_{g})^*, U) \cong {\rm
Ext}_{\Gamma}^{i+2}(K_{g+1}, U)=0$. Moreover, the exact sequence
$0 \to K_{g-1}\to ({\rm Im}f_{g})^* \to {\rm Ext}_{\Lambda}^{g}(M,
U) \to 0$ yields an exact sequence:
$${\rm Ext}_{\Gamma}^{g-1}(({\rm Im}f_{g})^*, U)
\to {\rm Ext}_{\Gamma}^{g-1}(K_{g-1}, U) \to {\rm
Ext}_{\Gamma}^{g}({\rm Ext}_{\Lambda}^{g}(M, U), U) \to {\rm
Ext}_{\Gamma}^{g}(({\rm Im}f_{g})^*, U).$$ So ${\rm
Ext}_{\Gamma}^{g}({\rm Ext}_{\Lambda}^{g}(M, U), U) \cong {\rm
Ext}_{\Gamma}^{g-1}(K_{g-1}, U) \cong {\rm
Ext}_{\Gamma}^{1}(K_{1}, U)$ and we get the following commutative
diagram with exact rows:

$$\xymatrix{P_1 \ar[d]^{\sigma _{P_{1}}} \ar[r]^{f_1} &
P_0 \ar[d]^{\sigma _{P_{0}}} \ar[r] & M \ar[d]^{h} \ar[r] & 0\\
P_{1}^{**} \ar[r]^{f_1^{**}} & P_{0}^{**} \ar[r] & {\rm
Ext}_{\Gamma}^{g}({\rm Ext}_{\Lambda}^{g}(M, U), U) \ar[r] & 0 }
$$
where $\sigma _{P_{1}}$ and $\sigma _{P_{0}}$ are isomorphisms.
Hence $h$ is also an isomorphism and ${\rm Ext}_{\Gamma}^{g}({\rm
Ext}_{\Lambda}^{g}(M, U),$
\linebreak
$U) \cong M(\neq 0)$. By
assumption, $_{\Lambda}U _{\Gamma}$ is a Gorenstein bimodule, so
grade$_{U}{\rm Ext}_{\Lambda}^{g}(M, U) \geq g$ and hence
grade$_{U}{\rm Ext}_{\Lambda}^{g}(M, U)=g$.

(2) It follows from (1) and its dual statement. $\blacksquare$

\vspace{0.2cm}

The following two corollaries are immediate consequence of Theorem
4.1.

\vspace{0.2cm}

{\bf Corollary 4.2} ([9] Proposition 3.1) {\it Let} $\Lambda$ {\it
and} $\Gamma$ {\it be left and right noetherian rings. If} $U
_{\Gamma}$ {\it is injective, then} $M\cong M^{**}$ {\it for any}
$M$ {\it in} mod $\Lambda$.

\vspace{0.2cm}

{\bf Corollary 4.3} {\it Under the assumptions of Theorem 4.1(2),
there is a duality between} $\mathcal{H}_{g}({\rm mod}\ \Lambda)$
{\it and} $\mathcal{H}_{g}({\rm mod}\ \Gamma ^{op})$ {\it given
by} $M \to {\rm Ext}_{\Lambda}^{g}(M, U)$ ({\it where} $M \in
\mathcal{H}_{g}({\rm mod}\ \Lambda)$).

\vspace{0.2cm}

The following result is a generalization of [16] Theorem 6, which
gives some characterizations of the modules in
$\mathcal{H}_{g}({\rm mod}\ \Lambda)$.

\vspace{0.2cm}

{\bf Theorem 4.4} {\it Let} $_{\Lambda}U _{\Gamma}$ {\it be a
Gorenstein bimodule with r.}id$_{\Gamma}(U)=${\it
l.}id$_{\Lambda}(U)=g$. {\it Then, for any} $0\neq M \in$ mod
$\Lambda$, {\it the following statements are equivalent.}

(1) $M \in \mathcal{H}_{g}({\rm mod}\ \Lambda)$.

(2) $M\cong {\rm Ext}_{\Gamma}^{g}({\rm Ext}_{\Lambda}^{g}(M, U),
U)$.

(3) $M\cong {\rm Ext}_{\Gamma}^{g}(N, U)$ {\it for some} $N\in$
mod $\Gamma ^{op}$.

(4) Hom$_{\Lambda}(M, \bigoplus _{i=0}^{g-1}E_{i})=0$.

\vspace{0.2cm}

{\it Proof.} $(1) \Rightarrow (2)$ follows from Corollary 4.3, and
$(2) \Rightarrow (3)$ is trivial.

$(3) \Rightarrow (4)$ Since $_{\Lambda}U _{\Gamma}$ is a
Gorenstein bimodule, $U$-resol.dim$_{\Lambda}(E_{i}) \leq i$ for
any $0 \leq i \leq g-1$ by Lemma 2.6. Then we get our conclusion
by Lemma 2.3 and Lemma 2.1.

$(4) \Rightarrow (1)$ Since {\it l.}id$_{\Lambda}(U)=g$,
Ext$_{\Lambda}^{i}(M, U)=0$ for any $i\geq g+1$. On the other
hand, we know that $M^*=0$ because Hom$_{\Lambda}(M, E_{0})=0$. In
addition, we have exact sequences:
$$0 \to K_{i-1} \to E_{i-1} \to K_{i} \to 0$$
for any $1\leq i \leq g-1$, where $K_{i-1}={\rm Ker}(E_{i-1}\to
E_{i})$. From Hom$_{\Lambda}(M, E_{i})=0$ we know that
Hom$_{\Lambda}(M, K_{i})=0$. But Hom$_{\Lambda}(M, K_{i})\to {\rm
Ext}_{\Lambda}^{i}(M, U)(\cong {\rm Ext}_{\Lambda}^{1}(M,
K_{i-1}))\to 0$ is exact, so ${\rm Ext}_{\Lambda}^{i}(M, U)=0$ for
any $1\leq i \leq g-1$ and grade$_{U}M\geq g$. We claim that
grade$_{U}M=g$. Otherwise, if grade$_{U}M>g$, we then have that
${\rm Ext}_{\Lambda}^{g}(M, U)=0$ and ${\rm Ext}_{\Lambda}^{i}(M,
U)=0$ for any $i\geq 1$. It follows from [15] Corollary 2.5 that
$M=0$, which is a contradiction. $\blacksquare$

\vspace{0.2cm}

Let $\mathcal{A}$ be an abelian category and $\mathcal{B}$ a full
subcategory of $\mathcal{A}$. An object $X \in \mathcal{A}$ is
called an embedding cogenerator for $\mathcal{B}$ if every object
in $\mathcal{B}$ admits an injection to some direct product of
copies of $X$ in $\mathcal{A} ^{[15]}$. For any $M$ in mod
$\Lambda$ we use $E(M)$ to denote the injective envelope of $M$.

\vspace{0.2cm}

{\bf Corollary 4.5} {\it Under the assumptions of Theorem 4.4,}
$E_{g}$ {\it is an injective embedding cogenerator for}
$\mathcal{H}_{g}({\rm mod}\ \Lambda)$.

\vspace{0.2cm}

{\it Proof.} Let $M$ be in $\mathcal{H}_{g}({\rm mod}\ \Lambda)$.
Notice that $M$ is finitely cogenerated, so, by [1] Proposition
18.18, $E(M) \cong E(S_{1})\bigoplus \cdots \bigoplus E(S_{t})$,
where each $S_{i}$ is isomorphic to a simple submodule of $M$ for
any $1 \leq i \leq t$. Since $M \in \mathcal{H}_{g}({\rm mod}\
\Lambda)$, each $S_{i} \in \mathcal{H}_{g}({\rm mod}\ \Lambda)$ by
Theorem 4.4.

Because {\it r.}id$_{\Gamma}(U)=g$, $\bigoplus _{i=0}^{g}E_{i}$ is
an injective embedding cogenerator for mod $\Lambda$ by [15]
Proposition 2.8. So Hom$_{\Lambda}(S_{i}, \bigoplus
_{i=0}^{g}E_{i})\neq 0$ and hence Hom$_{\Lambda}(S_{i}, E_{g})\neq
0$ by Theorem 4.4, which implies that each $S_{i}$ can be embedded
into $E_{g}$. Therefore $M \hookrightarrow E(M)\cong \bigoplus
_{i=0}^{t}E(S_{i})\hookrightarrow E_{g}^{(t)}$ and $E_{g}$ is an
injective embedding cogenerator for $\mathcal{H}_{g}({\rm mod}\
\Lambda)$. $\blacksquare$

\vspace{0.5cm}

\centerline{\bf 5. Finite injective dimension}

\vspace{0.2cm}

In this section we discuss the properties of $_{\Lambda}U
_{\Gamma}$ with finite left or right injective dimension. We first
have the following

\vspace{0.2cm}

{\bf Proposition 5.1} {\it If} {\it l.}id$_{\Lambda}(U)=k$ {\it
and} $E_{k}$ {\it is in} add$_{\Lambda}U$ ({\it equivalently},
$^*E_{k}$ {\it is flat}), {\it then} $_{\Lambda}U$ {\it is
injective}.

\vspace{0.2cm}

{\it Proof.} Assume that {\it l.}id$_{\Lambda}(U)=k\neq 0$. Then
there is a simple $\Lambda$-module $S$ such that
Ext$_{\Lambda}^{k}(S, U)\neq 0$. It is easy to see that
Hom$_{\Lambda}(S, E_{k})\cong$Ext$_{\Lambda}^{k}(S, U)$, so
Hom$_{\Lambda}(S, E_{k})\neq 0$ and hence there is an exact
sequence $0 \to S \buildrel {f} \over \longrightarrow E_{k} \to
{\rm Coker}f \to 0$, which yields an exact sequence ${\rm
Ext}_{\Lambda}^{k}(E_{k}, U)\to {\rm Ext}_{\Lambda}^{k}(S, U)\to
{\rm Ext}_{\Lambda}^{k+1}({\rm Coker}f, U)$. Since $E_{k} \in$
add$_{\Lambda}U$, ${\rm Ext}_{\Lambda}^{k}(E_{k}, U)=0$. On the
other hand, {\it l.}id$_{\Lambda}(U)=k$, so ${\rm
Ext}_{\Lambda}^{k+1}({\rm Coker}f, U)=0$. Hence ${\rm
Ext}_{\Lambda}^{k}(S, U)=0$, which is a contradiction.
$\blacksquare$

\vspace{0.2cm}

{\bf Corollary 5.2} {\it If l.}id$_{\Lambda}(\Lambda)=k$ {\it and
the} $(k+1)$st {\it term} ({\it that is, the last term}) {\it in a
minimal injective resolution of} $_{\Lambda}\Lambda$ {\it is flat,
then} $\Lambda$ {\it is self-injective}.

\vspace{0.2cm}

{\bf Corollary 5.3} {\it If l.}id$_{\Lambda}(U)=k<U$-${\rm
dom.dim}(_{\Lambda}U)$, {\it then} $_{\Lambda}U$ {\it is
injective}.

\vspace{0.2cm}

{\bf Corollary 5.4} ([17] Proposition 8) {\it If
l.}id$_{\Lambda}(\Lambda)=k<\Lambda$-${\rm
dom.dim}(_{\Lambda}\Lambda)$, {\it then} $\Lambda$ {\it is
self-injective}.

\vspace{0.2cm}

{\bf Proposition 5.5} {\it If} $_{\Lambda}U_{\Gamma}$ {\it is}
$k$-{\it Gorenstein and} {\it r.}id$_{\Gamma}(U)=${\it
l.}id$_{\Lambda}(U)=k$, {\it then}
$U$-resol.dim$_{\Lambda}(E_{k})$\linebreak
$=U$-resol.dim$_{\Gamma}(E'_{k})=k$ {\it and}
$_{\Lambda}U_{\Gamma}$ {\it is Gorenstein}.

\vspace{0.2cm}

{\it Proof.} Assume that $_{\Lambda}U$ is $k$-Gorenstein. By Lemma
2.6 and Lemma 2.3, {\it l.}fd$_{\Gamma}(^*E_i)\leq i$ for any $0
\leq i \leq k-1$. Since {\it r.}id$_{\Gamma}(U)=k$, there is a
module $X$ in mod $\Gamma ^{op}$ such that Ext$_{\Gamma}^{k}(X,
U)\neq 0$. Since $\bigoplus _{i=0}^{k}E_{i}$ is an injective
embedding cogenerator for mod $\Lambda$ by [15] Proposition 2.8,
it then follows from [6] Chapter VI, Proposition 5.3 that $0\neq
{\rm Hom}_{\Lambda}({\rm Ext}_{\Gamma}^k(X, U), \bigoplus _{i=0}^k
E_i)\cong {\rm Tor}_k^{\Gamma}(X, {^*(\bigoplus
_{i=0}^{k}E_i)})\cong \bigoplus _{i=0}^k{\rm Tor}_k^{\Gamma}(X,
{^*E_i})\cong {\rm Tor}_{k}^{\Gamma}(X, {^*E_k})$. So {\it
l.}fd$_{\Gamma}(^*E_k) \geq k$. On the other hand, by [14] Lemma
2.2, we have {\it r.}id$_{\Gamma}(U)=$sup$\{${\it
l.}fd$_{\Gamma}(^*E)| _{\Lambda}E$ is injective$\}$, so {\it
l.}fd$_{\Gamma}(^*E_k) \leq k$ and hence {\it
l.}fd$_{\Gamma}(^*E_k)=k$. By Lemma 2.3,
$U$-resol.dim$_{\Lambda}(E_{k})=k$. It then follows from Theorem
3.7 that $_{\Lambda}U$ is $(k+1)$-Gorenstein. In addition, {\it
l.}id$_{\Lambda}(U)=k$ by assumption, so $_{\Lambda}U$ is
Gorenstein. Similarly, we have that
$U$-resol.dim$_{\Gamma}(E'_{k})=k$ and $U_{\Gamma}$ is Gorenstein.
$\blacksquare$

\vspace{0.2cm}

Recall that an artin algebra is called an Auslander algebra if it
is $k$-Gorenstein for all $k$. By Proposition 5.5, we immediately
have the following

\vspace{0.2cm}

{\bf Corollary 5.6} ([8] Proposition 1.1) {\it If} $\Lambda$ {\it
is a} $k$-{\it Gorenstein algebra with right and left
selfinjective dimensions $k$, then the flat dimension of the}
$(k+1)$st {\it term in a minimal injective resolution of}
$_{\Lambda}\Lambda$ ({\it resp.} $\Lambda _{\Lambda}$) {\it is
equal to} $k$ {\it and} $\Lambda$ {\it is an Auslander algebra}.

\vspace{0.2cm}

Compare Corollary 5.3 with the following

\vspace{0.2cm}

{\bf Proposition 5.7} {\it If} {\it l.}id$_{\Lambda}(U)=k \leq
U$-${\rm dom.dim}(_{\Lambda}U)$, {\it then} $\bigoplus
_{i=0}^{k}E_{i}$ {\it is an injective embedding cogenerator for}
mod $\Lambda$ {\it if and only if} {\it r.}id$_{\Gamma}(U)=k$.

\vspace{0.2cm}

{\it Proof.} The sufficiency follows from [15] Proposition 2.8.
Now we prove the necessity. Since $U$-dom.dim$(_{\Lambda}U)\geq
k$, $E_{i}\in$add$_{\Lambda}U$ for any $0\leq i \leq k-1$. On the
other hand, {\it l.}id$_{\Lambda}(U)=k$ implies that $E_{i}=0$ for
any $i\geq k+1$. So $U$-resol.dim$_{\Lambda}(E_{k})\leq k$. Then,
by Lemma 2.3, {\it l.}fd$_{\Gamma}(^*E_k)\leq k$ and {\it
l.}fd$_{\Gamma}(^*E_i)=0$ for any $0\leq i \leq k-1$. It follows
that {\it l.}fd$_{\Gamma}(^*(\bigoplus _{i=0}^{k}E_i))\leq k$. So,
by Lemma 2.1, we have ${\rm Hom}_{\Lambda}({\rm
Ext}_{\Gamma}^{k+1}(X, U), \bigoplus _{i=0}^{k}E_{i})=0$ for any
$X\in$mod $\Gamma ^{op}$. However, $\bigoplus _{i=0}^{k}E_{i}$ is
an injective embedding cogenerator for mod $\Lambda$, so ${\rm
Ext}_{\Gamma}^{k+1}(X, U)=0$ and {\it r.}id$_{\Gamma}(U)\leq k$.
Hence we conclude that {\it r.}id$_{\Gamma}(U)=k$ by [4] Lemma
1.7. $\blacksquare$

\vspace{0.2cm}

Finally we conjecture the following, which is a generalization of
the Auslander-Reiten conjecture mentioned in Section 1: a
Gorenstein bimodule $_{\Lambda}U _{\Gamma}$ is cotilting, that is,
{\it l.}id$_{\Lambda}(U)<\infty$ and {\it
r.}id$_{\Gamma}(U)<\infty$.

\vspace{0.5cm}

{\bf Acknowledgements} The research of the author was partially
supported by Specialized Research Fund for the Doctoral Program of
Higher Education (Grant No. 20030284033). The author thanks the
referee for the useful comments.


\vspace{0.5cm}

\begin{thebibliography}{101}

\bibitem[1]{A1} F. W. Anderson and K. R. Fuller, Rings and Categories of
modules, 2nd ed, Graduate Texts in Mathematics {\bf 13},
Springer-Verlag, Berlin-Heidelberg-New York, 1992.

\bibitem[2]{A2} M. Auslander and M. Bridger, Stable Module Theory,
Memoirs Amer. Math. Soc. {\bf 94}, American Mathematical Society,
Providence, Rhode Island, 1969.

\bibitem[3]{A3} M. Auslander and I. Reiten, {\it Applications of
contravariantly finite subcategories}, Advances in Math. {\bf
86}(1991), 111--152.

\bibitem[4]{A4} M. Auslander and I. Reiten,
{\it Cohen-Macaulay and Gorenstein artin algebras}, in: G.O.
Michler and C.M. Ringel, eds. Representation Theory of Finite
Groups and Finite Dimensional Algebras, Bielefeld, 1991. Progress
in Mathematics {\bf 95}, Birkhauser, Basel, 1991, pp.221--245.

\bibitem[5]{A5} M. Auslander and I. Reiten,
{\it $k$-Gorenstein algebras and syzygy modules}, J. Pure and
Appl. Algebra {\bf 92}(1994), 1--27.

\bibitem[6]{A6} H. Cartan and S. Eilenberg, Homological Algebra,
Princeton University Press, Princeton, 1956.

\bibitem[7]{A7} R. M. Fossum, P.A. Griffith and I. Reiten, Trivial
Extensions of Abelian Categories, Lecture Notes in Mathematics
{\bf 456}, Springer-Verlag, Berlin-Heidelberg-New York, 1975.

\bibitem[8]{A8} K. R. Fuller and Y. Iwanaga, {\it On n-Gorenstein rings
and Auslander rings of low injective dimension}, Representations
of Algebras (Ottawa, ON, 1992), Canad. Math. Soc. Conf. Proc. {\bf
14}, Amer. Math. Soc., Providence, RI, 1993, pp. 175--183.

\bibitem[9]{A9} Z. Y. Huang, {\it On a generalization of
the Auslander-Bridger transpose}, Comm. Algebra {\bf 27}(1999),
5791--5812.

\bibitem[10]{A10} Z. Y. Huang, $\omega$-$k$-{\it torsionfree modules
and} $\omega$-{\it left approximation dimension}, Science in China
(Series A) {\bf 44}(2001), 184--192.

\bibitem[11]{A11} Z. Y. Huang, {\it Extension closure of relative
syzygy modules}, Science in China (Series A) {\bf 46}(2003),
611--620.

\bibitem[12]{A12} Z. Y. Huang, {\it Selforthogonal modules with
finite injective dimension II}, J. Algebra {\bf 264}(2003),
262--268.

\bibitem[13]{A13} Z. Y. Huang, {\it On} $U$-{\it dominant dimension},
J. Algebra {\bf 285}(2005), 669--681.

\bibitem[14]{A14} Z. Y. Huang, {\it Extension closure of relative}
$k$-{\it torsionfree modules}, Comm. Algebra (to appear).

\bibitem[15]{A15} Z. Y. Huang and G. H. Tang, {\it Self-orthogonal
modules over coherent rings}, J. Pure and Appl. Algebra {\bf
161}(2001), 167--176.

\bibitem[16]{A16} Y. Iwanaga, {\it Duality over Auslander-Gorenstein rings},
Math. Scand. {\bf 81}(1997), 184--190.

\bibitem[17]{A17} Y. Iwanaga and H. Sato, {\it Minimal injective
resolution of Gorenstein rings}, Comm. Algebra {\bf 18}(1990),
3835--3856.

\bibitem[18]{A18} T. Kato, {\it Rings of $U$-dominant
dimension}$\geq 1$, T\^{ohoku} Math. J. {\bf 21}(1969), 321--327.

\bibitem[19]{A19} T. Wakamatsu, {\it Tilting modules and Auslander's Gorenstein
property}, J. Algebra {\bf 275}(2004), 3--39.

\bibitem[20]{A20} R. Wisbauer, {\it Decomposition properties in
module categories}, Acta Univ. Carolinae--Mathematica et Physica
{\bf 26}(1985), 57--68.

\end{thebibliography}
\end{document}